\documentclass[leqno,12pt]{article}
\usepackage{amssymb,amsfonts}
\usepackage{amsmath,latexsym}
\usepackage{french}
\usepackage{amstext,epic,eepic,epsf,pslatex}
\usepackage{graphicx}

\newcommand{\bepr}{{\em Proof} } 
\newcommand{\enpr}{\hfill \rule{.5em}{.5em}}

\newcommand{\R}{{\mathbb R}}

\newcommand{\Tr}{\hbox{Tr\,}}

\def\Xint#1{\mathchoice
{\XXint\displaystyle\textstyle{#1}}%
{\XXint\textstyle\scriptstyle{#1}}%
{\XXint\scriptstyle\scriptscriptstyle{#1}}%
{\XXint\scriptscriptstyle\scriptscriptstyle{#1}}%
\!\int}
\def\XXint#1#2#3{{\setbox0=\hbox{$#1{#2#3}{\int}$ }
\vcenter{\hbox{$#2#3$ }}\kern-.6\wd0}}

\newtheorem{defin}{Definition}[section] 
\newtheorem{prop}{Proposition}[section] 
\newtheorem{thm}{Theorem}[section] 
\newtheorem{lemma}{Lemma}[section]

\begin{document}

\title{Compensated integrability.  Applications to the Vlasov--Poisson equation \\ 
and other models in mathematical physics}

\author{Denis Serre \\ \'Ecole Normale Sup\'erieure de Lyon\thanks{U.M.P.A., UMR CNRS--ENSL \# 5669. 46 all\'ee d'Italie, 69364 Lyon cedex 07. France. {\tt denis.serre@ens-lyon.fr}}\thanks{The author thanks Meiji University (Tokyo) and the Chinese University (Hong Kong) for their hospitality when a part of this paper was developed.}}

\date{}

\maketitle

\begin{english}
\begin{abstract}
We extend our analysis of divergence-free positive symmetric tensors (DPT) begun in a previous paper. On the one hand, we refine the statements and give more direct proofs. Next, we study the most singular DPTs, and use them to prove that the determinant is the only quantity that enjoys an improved integrability. Curiously, these singularities are intimately related to the Minkowski's Problem for convex bodys with prescribed Gaussian curvature. We then cover a list of models of mathematical physics that display a divergence-free symmetric tensor~; the most interesting one is probably that of nonlinear Maxwell's equations in a relativistic frame. The case of the wave equation is the occasion to highlight the role of the positivity assumption. Last, but not least, we show that the Vlasov--Poisson equation for a plasma is eligible for our theory.
\end{abstract}
\end{english}

\begin{resume}
Nous poursuivons l'analyse des tenseurs sym\'etriques positifs \`a divergence nulle (DPT) commenc\'ee dans un pr\'ec\'edent article. Pour une part, nous \'etablissons des \'enonc\'es plus fins, et les preuves sont plus directes. Dans un second temps, nous \'etudions les DPTs les plus singuliers, et nous les utilisons pour montrer que le d\'eterminant est la seule quantit\'e qui jouisse d'une int\'egrabilit\'e plus \'elev\'ee. Curieusement, ces singularit\'es sont intimement li\'ees au Probl\`eme de Minkowski, qui concerne les corps convexes \`a courbure de Gau\ss\, prescrite comme fonction de la direction normale. Nous passons ensuite en revue un certain nombre de mod\`eles de Physique Math\'ematique o\`u les tenseurs sym\'etriques \`a divergence nulle sont \`a l'{\oe}uvre. Le plus int\'eressant d'entre eux est sans doute le syst\`eme de Maxwell dans un contexte non-lin\'eaire et relativiste. Le cas de l'\'equation des ondes montre que l'hypoth\`ese de positivit\'e est essentielle pour la th\'eorie. Enfin, nous montrons comment le cadre des DPTs s'applique \`a l'\'equation de Vlasov--Poisson avec force r\'epulsive (cas d'un plasma).
\end{resume}



\section{Introduction}

Let $d\ge2$ be an integer. In a previous paper \cite{Ser_DPT}, we began the study of divergence-free positive symmetric tensors, in short DPTs. These objects are measurable fields
$$T:\Omega\longrightarrow{\bf Sym}_d^+$$
over an open domain $\Omega\subset\R^d$ (or a torus $\Omega=\R^d/\Gamma$), with values in the cone of positive semi-definite matrices, whose row-wise divergence ${\rm Div}\,T$ vanishes:
$$\forall\,i=1,\ldots,d\qquad\sum_{j=1}^d\partial_jT_{ij}\equiv0.$$
We also considered situations where the condition above is replaced by some control of ${\rm Div}\,T$ in the space of vector fields with values in the space of bounded measures.

In the periodic case, our main result was that when such a field is  integrable, then $(\det T)^{\frac1{d-1}}$ is integrable too. This is remarkable since a product like $(T_{11}\cdots T_{dd})^{\frac1{d-1}}$ does not need to be integrable. Such a result can therefore be called a {\em Compensated Integrability}. At the algebraic level, the relation between the differential constraint and the function $\det^{\frac1{d-1}}$ is that the latter is concave in the characteristic directions of the operator ${\rm Div}$. The inequality associated with the property of compensated integrability,
\begin{equation}
\label{eq:funcin}
\Xint-_{\R^d/\Gamma}(\det T)^{\frac1{d-1}}dx\le\left(\det\Xint-_{\R^d/\Gamma}T(x)\,dx\right)^{\frac1{d-1}},
\end{equation}
is a {\em non-diagonal} extension of Gagliardo's Inequality.
It tells us that the map  $T\longmapsto(\det T)^{\frac1{d-1}}$ is {\em Divergence-quasiconcave}, in the terminology of Fonseca \& M\"uller \cite{FM}.

For a bounded domain, we established an inequality similar to (\ref{eq:funcin}), where the right-hand side involves instead a boundary integral. Despite its correctness, our Theorem 2.3 was still perfectible. On the one hand, we proved only the functional inequality, but felt short of establishing the gain of integrability and had to make it an assumption. On the other hand, we imposed the useless restriction that the domain $\Omega$ be convex.

\bigskip

The paper below addresses several aspects of DPTs, concerning either general statements or their applications to Mathematical Physics and Fluid Dynamics.

In Section \ref{s:CI}, we complete the results of \cite{Ser_DPT} and give more elegant proofs, both in the periodic and in the bounded cases. We drop the assumption of convexity and establish the Compensated Integrability in all situations. We treat the case where $\Omega$ is a slab $(0,\tau)\times\R^{d-1}$, which is taylored for the applications to models of continuum mechanics, especially inviscid compressible fluids~; our statement can be applied whenever the total mass and energy of the fluid are finite. As shown in \cite{Ser_DPT}, this yields new {\em a priori} estimates, where the gain of integrability comes to the price of a time integration. We address a question that was left aside,  of whether the positivity of $T(x)$ is a meaningfull assumption. The answer is positive when $d\ge3$, as shown by an example, constructed by S. Klainerman \& M. Machedon in the context of the wave equation.

\bigskip

Section \ref{s:Imm} is a rather general study of those functions $T\longmapsto f(T)$ which display a similar gain of integrability as $\det^{\frac1{d-1}}$. We first show that they must vanish on rank-one tensors $v\otimes v$. It is therefore natural to consider appropriate powers of {\em immanants}. Our main result is that among this collection, $\det^{\frac1{d-1}}$ is the only one to gain maximal integrability. This paragraph makes use of the classification of the locally integrable DPTs that are homogeneous of a given degree. On the one hand, we show that this degree may not belong to the interval $(-d,1-d)$~; this is an indication that singularities of DPTs should not be too strong. On the other hand, we characterize those of degree $1-d$, which are intimately related to a geometrical problem posed by Minkowski.

\bigskip

Section \ref{s:otherDPT} is devoted to divergence-free symmetric tensors that occur in Mathematical Physics~; they often come as energy-momentum tensors. We warn the reader that not all of them are positive, and thus we do not always expect a compensated integrability. We mentionned above the case of the wave equation. We  present here a few other examples, among which a rather interesting one is the non-linear Maxwell system of an electro-magnetic field~; it is remarkable that the symmetry of the associated tensor is a consequence of the invariance of the model under the action of the Lorentz group of special relativity. We  also describe in full generality those first-order systems of conservation laws having the DPT form, and that are compatible with a convex entropy~; the Euler system of an inviscid gas falls into this category. We show that such systems can be described in terms of a single potential function.

\bigskip

Last, but not least, Section \ref{s:Vla} deals with kinetic models of the Vlasov family in $d=1+n$ space-time dimensions. We show that such equations are associated with a divergence-free symmetric tensor, which turns out to be positive when the self-induced force $F$ is repulsive. This covers the case of a plasma (Coulomb force), and excludes that of a galaxy (gravity force). The construction is rather indirect in terms of the mass density, but this is the price to pay for the positivity of the tensor.  As in gas dynamics, this leads us to a new  estimate, which displays a higher integrability property in space and time.   A previous, more naive attempt in \cite{Ser_DPT}, led us to a tensor containing a diagonal  block $F\otimes F-\frac12\,|F|^2I_n$, which is neither positive nor negative, so that compensated integrability could not be applied.

\section{Compensated integrability}\label{s:CI}

\subsection{Periodic case}

Let us recall  Theorem 2.1 of \cite{Ser_DPT}~:
\begin{thm}\label{th:deun}
Let the DPT $x\longmapsto A(x)$ be $\Gamma$-periodic, with $A\in L^1(\R^d/\Gamma)$. Then $(\det A)^{\frac1{d-1}}\in L^1(\R^d/\Gamma)$ and there holds
\begin{equation}\label{Jdetper}
\Xint-_{\R^d/\Gamma}(\det A(x))^{\frac1{d-1}}dx\le\left(\det\Xint-_{\R^d/\Gamma}A(x)\,dx\right)^{\frac1{d-1}}.
\end{equation}
\end{thm}

The fact that $(\det A)^{\frac1{d-1}}$ is integrable, is not implied by the sole assumption $A\in L^1(\R^d/\Gamma)$. There is no reason why other homogeneous functions of $A$ of degree $\frac d{d-1}$ should be integrable. The theorem therefore expresses a compensation property, which  may be called {\em Compensated Integrability}. 

A similar phenomenon was established by Coifman \& al. \cite{CLMS}, who proved that if two vector fields $u,v\in L^2(\R^d)$ satisfy ${\rm div}\,u=0$ and ${\rm curl}\,v=0$, then the scalar product $u\cdot v$ belongs to the Hardy space ${\cal H}^1(\R^d)$, a strict subspace of the obvious $L^1(\R^d)$. If $f=u\cdot v$ is non-negative, this implies that 
$f\log(1+f)\in L^1_{\rm loc}(\R^d)$. It is interesting to compare the result of Coifman \& al. with ours when $d=2$, because their differential constraints read ${\rm Div}\,T=0$ for $T=\begin{pmatrix} u_1 & u_2 \\ -v_2 & v_1 \end{pmatrix}$, the function $f=u\cdot v$ is nothing but $\det T$, and we both make the assumption that $f\ge0$. Their integrability is higher than ours, both in the assumption and in the conclusion. We have an extra assumption (symmetry), but our gain is better (from $L^{1/2}$ to $L^1$, instead of from $L^1$ to $L\log L$).

\bigskip

The following proof of Theorem \ref{th:deun} is much more direct that the one in \cite{Ser_DPT}. Let $f>0$ be an arbitrary smooth, $\Gamma$-periodic function. Let $S\in{\bf SPD}_d$ satisfy the constraint
\begin{equation}
\label{eq:detS}
\det S=\Xint-_{\R^d/\Gamma}f(x)\,dx.
\end{equation}
By Yan Yan Li's Theorem \cite{YYL}, the Monge--Amp\`ere equation $\det {\rm D}^2\theta=f$ admits a smooth, convex solution of the form $\theta=\frac12x^TSx+\rho(x)$ where $\rho$ is periodic too. With the same arguments as in  \cite{Ser_DPT}, we have
$$(f\det A)^{\frac1d}\le\frac1d\,\Tr(A{\rm D}^2\theta)=\frac1d\,(\Tr(AS)+{\rm div}(A\nabla\rho)).$$
Integrating, we have
$$\Xint-_{\R^d/\Gamma}(f\det A)^{\frac1d}dx\le\frac1d\,\Tr(S\Xint-_{\R^d/\Gamma}A(x)\,dx)=:\frac1d\,\Tr(SA_+).$$
We optimize the choice of $S$ by taking $S=\lambda A_+^{-1}$, where $\lambda>0$ satisfies (\ref{eq:detS}). We obtain
$$\Xint-_{\R^d/\Gamma}(f\det A)^{\frac1d}dx\le\lambda=\left(\Xint-_{\R^d/\Gamma}f(x)\,dx\,\det A_+\right)^{\frac1d}.$$
Denoting $\phi:=f^{\frac1d}$, this is recast as
$$\Xint-_{\R^d/\Gamma}\phi(\det A)^{\frac1d}dx\le\lambda=\|\phi\|_{L^d_{\rm per}}\,(\det A_+)^{\frac1d}.$$
By density, this inequality remains true for every non-negative function $\phi\in L^p_{\rm per}$. Because $\det A\ge0$, we infer that $(\det A)^{\frac1d}\in L^{d'}_{\rm per}$, where $d'=\frac d{d-1}$ is the conjugate exponent, that is $(\det A)^{\frac1{d-1}}\in L^1_{\rm per}$, and we have
$$\|(\det A)^{\frac1d}\|_{L^{d'}_{\rm per}}\le(\det A_+)^{\frac1d}.$$
This is exactly our functional inequality (\ref{Jdetper}).

\bigskip

\paragraph{Remark.} The proof given here has the advantage to apply under a weaker hypothesis, when we assume that the entries of $A$ are bounded measures on $\R^d/\Gamma$, instead of being integrable. The positiveness means that for every $\xi\in S^{d-1}$, the bounded measure $\sum_{i,j}\xi_i\xi_ja_{ij}$ is non-negative. The only point to clarify is the definition of $(\det A)^{\frac1d}$. By assumption, the entries are absolutely continuous with respect to the non-negative measure $\Tr A$. One can therefore introduce the densities $f_{ij}$ by $a_{ij}=f_{ij}\,\Tr A$. These functions $f_{ij}$ are integrable with respect to $\Tr A$ and the matrix $F=(f_{ij})_{1\le i,j\le d}$ takes values in ${\bf Sym}_d^+$. Then 
$$(\det A)^{\frac1d}:=(\det F)^{\frac1d}\Tr A$$
defines a bounded measure.
Once again, there is a compensated integrability: $(\det A)^{\frac1d}$ is absolutely continuous with respect to the Lebesgue measure, and its density belongs to $L^{\frac d{d-1}}(\R^d/\Gamma)$.

The present remark can be made for the Theorems \ref{th:improv} and \ref{th:slab} below. It is also meaningful when treating renormalized solutions of Boltzman flow (see Section 3.4 of \cite{Ser_DPT}), where the conservation of momentum includes a {\em defect measure} $\Sigma$, since our analysis provided an estimate of 
$$\int_0^\tau\int_{\R^n}(\rho\det\Sigma)^{\frac1n}dy\,dt.$$

\subsection{DPTs in a bounded domain}

The analysis above works the same in a convex bounded domain $\Omega$. For the sake of completeness, we may assume that ${\rm Div}\,A$ is a vector-valued bounded measure. Let $f>0$ be an arbitrary smooth function over $\overline\Omega$. Let $\nabla\theta$ be the Brenier's transport from $(\Omega,f(x)dx)$ to $(B_r,dx)$, where the constraint
$$\frac{r^d}d|S^{d-1}|=\int_\Omega f(x)\,dx$$
expresses that both measures have the same mass.

Following the same lines, we have
$$(f\det A)^{\frac1d}\le\frac1d\,({\rm div}(A\nabla\theta)+({\rm Div}\,A)\nabla \theta).$$
Integrating, there comes
$$\int_\Omega(f\det A)^{\frac1d}dx\le\frac1d\,(\int_{\partial\Omega}(A\nabla\theta)\cdot\vec n\,ds(x)+\int_\Omega({\rm Div}\,A)\nabla \theta\,dx)\le\frac rd\,(\int_{\partial\Omega}|A\vec n|\,ds(x)+\|{\rm Div}\,A\|_{\cal M}).$$
Let us proceed as above: we introduce $\phi=f^{\frac1d}$ and rewrite
$$\int_\Omega\phi(\det A)^{\frac1d}dx\le K\,\|\phi\|_{L^d}$$
where
$$K=\frac1d\,\left(\frac d{|S^{d-1}|}\right)^{\frac1d}(\int_{\partial\Omega}|A\vec n|\,ds(x)+\|{\rm Div}\,A\|_{\cal M}).$$
We deduce that $(\det A)^{\frac1d}\in L^{d'}(\Omega)$, with $\|(\det A)^{\frac1d}\|_{L^{d'}}\le K$. This proves that $(\det A)^{\frac1{d-1}}$ is integrable, a point that was left open in Theorem 2.3 of \cite{Ser_DPT}.


\subsubsection{Non-convex domains}

This Theorem 2.3 can actually be improved in another direction, by dropping the assumption that the domain be convex. Our ultimate result is therefore:
\begin{thm}\label{th:improv}
Let $\Omega$ be a bounded open subset in $\R^d$ with smooth boundary. Let $A\in L^1(\Omega;{\bf Sym}_d^+)$ be given, such that ${\rm Div}\,A$ is a (vector-valued) bounded measure. We assume that the normal trace $A\vec n$ is a bounded measure. Then there holds
\begin{equation}\label{convmeas}
\int_{\Omega}(\det A(x))^{\frac1{d-1}}dx\le\frac1{d|S^{d-1}|^{\frac1{d-1}}}\,\left(\|A\vec n\|_{{\cal M}(\partial\Omega)}+\|{\rm Div}\,A\|_{{\cal M}(\Omega)}\right)^{\frac d{d-1}}.
\end{equation}
\end{thm}

\bepr

The Theorem has already been proven when the domain is convex and the normal trace is integrable. We may therefore apply this restricted version to the pair $(\tilde\Omega,\tilde A)$, where $\tilde\Omega$ is a ball that contains $\overline\Omega$, and $\tilde A$ is the extension of $A$ by $0_d$ over $\tilde\Omega\setminus\Omega$. The assumption tells us that ${\rm Div}\tilde\Omega$ is a bounded measure, and we have
$$\|{\rm Div}\,\tilde A\|_{{\cal M}(\tilde\Omega)}=\|A\vec n\|_{{\cal M}(\partial\Omega)}+\|{\rm Div}\,A\|_{{\cal M}(\Omega)}.$$
Besides, the normal trace of $\tilde A$ vanishes identically. The inequality (\ref{convmeas}), valid for $\tilde A$ in $\tilde\Omega$, yields the same inequality for $A$ in $\Omega$.

\enpr

\bigskip

Of course, the same remark as in the periodic case holds true: we can actually replace the integrability of $A$, by the assumption that $A$ is a symmetric tensor of bounded measures, non-negative in the sense that for every vector $\xi\in\R^d$, the scalar measure $\xi^T A\xi$ is non-negative. This assumption allows us to define a non-negative measure $(\det A)^{\frac1d}$, and the theorem tells us that the latter is actually a function of class $L^{\frac d{d-1}}$. In particular, it does not display a singular part.

\bigskip

Remark also that the isoperimetric inequality follows immediately and for every smooth domain, by applying Theorem \ref{th:improv} to the  tensor $A\equiv I_d$.

\subsection{The case of a slab}\label{s:slab}

We consider now a domain of the form $\Omega=(0,\tau)\times\R^n$, where $\tau>0$ and $d=1+n$. We split the coordinates as $x=(t,y)$ and may think of $t$ and $y$ as  time and space variables, respectively. Suppose that $A$ is a DPT over $\Omega$, with $A\in L^1(\Omega)$. Assume also that its normal traces on the top/bottom boundaries are bounded measures. 

Let us choose a function $\phi\in{\cal D}(\R^n)$ with support in the ball $B_2$ of radius $2$, which satisfies $0\le\phi\le1$ everywhere, and $\phi\equiv 1$ in the unit ball $B_1$. If $R>0$, define $\phi_R(y):=\phi(y/R)$.
Let us apply Theorem \ref{th:improv} to the tensor $\phi_RA$, in the domain $\Omega_R=(0,\tau)\times B_{2R}$. Its divergence is $A\nabla\phi_R$ and its trace on the lateral boundary $(0,\tau)\times\partial B_{2R}$ vanishes, while $|\phi_RA\xi|\le|A\xi|$. We have therefore
\begin{eqnarray*}
\int_{\Omega_R}(\phi_R\det A(x))^{\frac1{d-1}}dx & \le & \frac1{d|S^{d-1}|^{\frac1{d-1}}}\,\left(\|\phi_RA\vec n\|_{{\cal M}(\partial\Omega_R)}+\|{\rm Div}\,(\phi_RA)\|_{{\cal M}(\Omega_R)}\right)^{\frac d{d-1}} \\
& \le & \frac1{d|S^{d-1}|^{\frac1{d-1}}}\,\left(\|A\vec n\|_{{\cal M}(\{0\}\times\R^n)}+\|A\vec n\|_{{\cal M}(\{\tau\}\times\R^n)}+\|A\nabla\phi_R\|_{{\cal M}(\Omega_R)}\right)^{\frac d{d-1}}. 
\end{eqnarray*}
We infer
$$\int_0^\tau\int_{B_R}(\det A)^{\frac1{d-1}}dy\,dt\le \frac1{d|S^{d-1}|^{\frac1{d-1}}}\,\left(\|A\vec n\|_{{\cal M}(\{0\}\times\R^n)}+\|A\vec n\|_{{\cal M}(\{\tau\}\times\R^n)}+\frac cR\,\|A\|_{L^1}\right)^{\frac d{d-1}}.$$
Passing to the limit as $R\rightarrow+\infty$, we may state:
\begin{thm}\label{th:slab}
Let $A$ be a DPT over the slab $(0,\tau)\times\R^n$. We assume that it is integrable, and that its normal traces on the top/bottom boundaries are bounded measures. Then $(\det A)^{\frac1n}$ is integrable and one has
\begin{equation}
\label{eq:slab}
\int_0^\tau\int_{\R^n}(\det A)^{\frac1n}dy\,dt\le \frac1{(n+1)|S^n|^{\frac1n}}\,\left(\|A\vec n\|_{{\cal M}(\{0\}\times\R^n)}+\|A\vec n\|_{{\cal M}(\{\tau\}\times\R^n)}\right)^{1+\frac1n}.
\end{equation}
\end{thm}

\bigskip

This is the fundamental inequality that we apply to various models in physics, mostly in gas dynamics.

\subsection{More about the gain integrability}

Theorems \ref{th:deun}, \ref{th:improv} and \ref{th:slab} turn out to be understatements in some circumstances. When we say that $(\det A)^{\frac1{d-1}}$ is integrable, we really mean that the function $x\mapsto(\det A(x))^{\frac1d}$, which we already know to be integrable, has the additional property to belong to $L^{\frac d{d-1}}$ (see the proofs). Of course, this implies that the function $x\mapsto(\det A(x))^{\frac1{d-1}}$, which is defined almost everywhere and measurable, is integrable. But this is not the end of the story. We shall encounter in the next section singular (though integrable) DPTs for which it makes sense to consider a more accurate definition of $(\det A)^{\frac1{d-1}}$, which includes a singular measure. Typically, if $\theta$ is a convex, positively homogeneous function of order one (for instance a semi-norm), the tensor $A:=\widehat{{\rm D}^2\theta}$ is locally integrable, rank-one away from the origin and therefore satisfies $\det A\equiv0$ almost everywhere. This violates the  identity
$$\int_B(\det A)^{\frac1{d-1}}dx=\int_B\det{\rm D}^2\theta\,dx={\rm vol}(\nabla\theta(B)),$$
where $B$ is the unit ball, and $\nabla\theta(B)$ must be understood as a convex body (approach $\theta$ by smooth convex functions $\theta_\epsilon$, or consider $\nabla\theta$ as a sub-differential). This strongly suggests to define, in this example,
$$(\det A)^{\frac1{d-1}}={\rm vol}(\nabla\theta(B))\cdot\delta_{x=0}.$$
Such an extension of the definition of the now bounded measure $(\det A)^{\frac1{d-1}}$ is meaningful only when $A$ displays the highest singularities (see Proposition \ref{p:sing} below), typically when $A$ itself is not locally in $L^{\frac d{d-1}}$. It is consistent with the non-diagonal Gagliardo inequalities stated above. We intend to develop more on this subject in a future work.

\section{Singularities {\em vs} gain of integrability}\label{s:Imm}

This section intends to show that the determinant is essentially the only homogeneous polynomial over ${\bf Sym}_d^+$ to display a maximal gain of integrability when applied to DPTs.

\subsection{Homogeneous singularities of DPTs}

We begin by considering those DPTs that are homogeneous functions of $x$, with a negative  degree $-m$~; we call $m$ the {\em order of singularity}. Denoting $r=|x|$ and $e=\frac xr$\,, we look for tensors of the form
$$T(x)=r^{-m}S(e)$$
where $e\longmapsto S(e)$ takes values in ${\bf Sym}_d^+$.
Of course, we restrict our attention to integrable singularities, more generally to singularities that are bounded measures. We assume therefore that $m<d$ and $S$ is a bounded measure over the unit sphere $S^{d-1}$. 

We begin with elementary examples:
\begin{lemma}
When $m\le d-1$, the tensor 
$$T_m=\frac1{r^m}\,(m\,e\otimes e+(d-1-m)I_d)$$
is a locally integrable DPT.
\end{lemma}
The family $T_m$ includes the special cases
$$T_0=(d-1)I_d,\qquad T_{d-1}=\frac{d-1}{r^{d-1}}\,e\otimes e,$$
where the constant factors $d-1$ are irrelevant. 

\bigskip

\bepr

Using the identity ${\rm Div}\,(fS)=f\,{\rm Div}\,S+S\nabla f$ and the fact ${\rm Div}(x\otimes x)=(d+1)x$, we find that ${\rm Div}\,T_m\equiv0$ away from the origin. When $\vec\phi\in{\cal D}(\R^d)^d$ and $\epsilon\in(0,1)$, we have
\begin{eqnarray*}
\langle{\rm Div}\,T_m,\vec\phi\rangle & = & -\langle T_m,\nabla\vec\phi\rangle=-\int_{\R^d}T_m:\nabla\vec\phi\,dx=-\int_{B_\epsilon}T_m:\nabla\vec\phi\,dx+\int_{S_\epsilon}(T_me)\cdot\vec\phi\,ds(e) \\
& = & O(\epsilon^{d-1-m}).
\end{eqnarray*}
Letting $\epsilon\rightarrow0^+$, we infer $\langle{\rm Div}\,T_m,\vec\phi\rangle=0$, that is ${\rm Div}\,T_m=0$, in the distributional sense.

\enpr

\bigskip

Remark that if $m\in(d-1,d)$, the divergence-free tensor $T_m$ is integrable but not positive. This is consistent with the first part of the following statement.
\begin{prop}\label{p:sing}
\begin{enumerate}
\item
In the class of homogeneous DPTs, the order of singularity $m< d$ must  be less than or equal to $d-1$.
\item
The  homogeneous DPTs  with order of singularity  $d-1$ are of the form
\begin{equation}
\label{eq:dmsun}
T(x)=\frac{\lambda(e)}{r^{d-1}}\,e\otimes e
\end{equation}
for some non-negative measure $\lambda$ satisfying
\begin{equation}
\label{eq:intla}
\int_{S^{d-1}}e_id\lambda(e)=0,\qquad\forall\,i=1,\ldots,d.
\end{equation}
For every such measure, the formula above provides a DPT whose entries are bounded measures.

More generally, a tensor given by Formula (\ref{eq:dmsun}) satisfies
$${\rm Div}\,T=\delta_{x=0}\otimes V_\lambda,\qquad V_\lambda:=\int_{S^{d-1}}e\,d\lambda(e)$$
in the sense of distributions.
\item
If $\theta:\R^d\rightarrow\R$ is convex and positively homogeneous of degree one (for instance if $\theta$ is a semi-norm), then  $\widehat{{\rm D}^2\theta}$ is a homogeneous DPT of maximal order of singularity, that is $d-1$.
\item 
Conversely, let $T$ be given by (\ref{eq:dmsun}) where $\lambda$ is strictly positive and satisfies (\ref{eq:intla}). Suppose that $\lambda\in{\cal C}^m(S^{d-1})$ for some $m\ge3$. Then there exists a convex function $\theta\in{\mathcal C}^{m+1,\alpha}(\R^d)$ (for every $\alpha\in(0,1)$), positively homogeneous of degree one, such that $T=\widehat{{\rm D}^2\theta}$.
\end{enumerate}
\end{prop}

\bigskip

\paragraph{Remark.} In point $4$, the map $e\mapsto\nabla \theta(e)$ solves Minkowski's Problem, which consists in finding a convex hypersurface $\Sigma$ whose Gau\ss\, curvature, here $\kappa=\frac1{\lambda(e)}$\,, is a prescribed function of the normal direction. The function $\theta$ is the support function of the body whose boundary is $\Sigma$. The solution of  the problem is related to an equation of Monge-Amp\`ere type on the sphere, namely 
$$\det{\rm D}^2\theta|_{e^\bot}=\lambda(e),$$
where the convexity of $\theta$ as a function of $x\in\R^d$ ensures the ellipticity. If $\lambda$ is even, the condition (\ref{eq:intla}) is obviously satisfied and the solution $\theta$ is a semi-norm.

\bigskip

\bepr

Let $T$ be a homogeneous DPT of order $m<d$. The divergence-free condition writes
$${\rm Div}\,S=\frac m{r^2}\,Sx,$$
where we view $S$ as a function of $x$, homogeneous of degree zero. Let us make the scalar product with $x$,
$$me^TSe=x\cdot{\rm Div}\,S={\rm div}(Sx)-\Tr S.$$
Let us  integrate over the unit ball $B$,
$$\int_B(me^TSe+\Tr S)\,dx=\int_{S^{d-1}}e^TSe\,ds(e).$$
With $dx=r^{d-1}dr\,ds(e)$, this yields
$$\int_{S^{d-1}}((m+1-d)e^TSe+\Tr S|_{e^\bot})\,ds(e)=0.$$
When $m\ge d-1$, the integrand is non-negative (recall that $S(e)\in{\bf Sym}_d^+$) and the equality above is equivalent to $(m+1-d)e^TSe+\Tr S|_{e^\bot}\equiv0$. If $m>d-1$, this means $S\equiv0_d$, which proves the first item. If $m=d-1$, it says only that $S$ vanishes over $e^\bot$, that is $S$ is of the form $\lambda(e)e\otimes e$ with $\lambda(e)\ge0$.

Conversely, let $\lambda$ be a non-negative measure, so that $T:=r^{1-d}\lambda(e)e\otimes e$ is a symmetric, positive semi-definite tensor. If $\lambda\equiv1$, we verify easily that the corresponding tensor $T^0$ equals   $\widehat{{\rm D}^2\theta}$ where $\theta\equiv|x|$. Therefore $T^0$ is a DPT. For a general $\lambda$, we have
$${\rm Div}\,T={\rm Div}(\lambda T^0)=T^0\nabla\lambda=r^{-d}(x\cdot\nabla\lambda)e\equiv0$$
away from the origin, where the last identity is that of Euler for homogeneous functions. If $\vec\phi\in{\cal D}(\R^d)^d$ and $\epsilon\in(0,1)$, we have
\begin{eqnarray*}
\langle{\rm Div}\,T,\vec\phi\rangle & = & -\langle T,\nabla\vec\phi\rangle=-\int_{\R^d}T:\nabla\vec\phi\,dx=-\int_{B_\epsilon}T:\nabla\vec\phi\,dx+\int_{S_\epsilon}(Te)\cdot\vec\phi\,ds(e) \\
& = & O(\epsilon^{d-m})+\epsilon^{1-d}\int_{S^{d-1}}e\cdot\vec\phi\,d\lambda(e).
\end{eqnarray*}
Letting $\epsilon\rightarrow0^+$, we obtain
$$\langle{\rm Div}\,T,\vec\phi\rangle=\vec\phi(0)\cdot\int_{S^{d-1}}e\,d\lambda(e),$$
that is ${\rm Div}\,T=\delta_{x=0}\otimes V_\lambda$ in the sense of distributions.
Therefore $T$ is a DPT if and only if (\ref{eq:intla}) is satisfied.

Next, if $\theta$ is convex and positively homogeneous of degree one, then we already know that   $T=\widehat{{\rm D}^2\theta}$ is a DPT. Since  $\widehat{{\rm D}^2\theta}$ is homogeneous of degree $-1$, $T$ is homogeneous of degree $1-d$.

The last point is Theorem 1, page 33 of Pogorelov \cite{Pog}.

\enpr

\subsection{Are there other gains of integrability ?}

\begin{defin}
Let $f:{\bf Sym}_d^+\longrightarrow\R^+$ be a continuous, homogeneous function of degree $\frac d{d-1}$\,. We say that $f$ exhibits a {\em maximal gain of regularity} for DPTs is there exists a finite constant $c_f$ such that for every divergence-free $T\in L^1(B;{\bf Sym}_d^+)$ with an integrable normal trace, one has
\begin{equation}
\label{eq:gain}
\int_Bf(T)\,dx\le c_f\|T\vec n\|_{L^1(S^{d-1})}^{\frac d{d-1}}.
\end{equation}
\end{defin}

The degree of homogeneity is justified by the fact that if {\em all} the derivatives $\partial T$ are under control, and not only ${\rm Div}\,T$,  then $T\in W^{1,1}(B)$, which is contained in $L^{\frac d{d-1}}(B)$ by Sobolev embedding. Besides, a dimensional analysis shows that such an inequality would be valid in every ball of $\R^d$, with the same constant $c_f$.

\bigskip

We begin by testing (\ref{eq:gain}) against the tensors $T_m$. Because of  $\vec n=e$, we have $T_m\vec n=(d-1)e$ and thus $\|T_m\vec n\|_{L^1(S^{d-1}}=(d-1)|S^{d-1}|$. On the other hand,
$$\int_Bf(T_m)\,dx=\frac{d-1}{d(d-1-m)}\,\int_{S^{d-1}}f(m\,e\otimes e+(d-1-m)I_d)\,ds(e).$$
We deduce therefore
\begin{lemma}\label{l:gTm}
If $f$ displays a maximal gain of integrability, then there exists a finite constant $c_f'$ such that
\begin{equation}
\label{eq:gTm}
\int_{S^{d-1}}f(m\,e\otimes e+(d-1-m)I_d)\,ds(e)\le c_f'(d-1-m)
\end{equation}
for every $m< d-1$.

By letting $m\rightarrow d-1$, $f$ must vanish on the cone of rank-one tensors $v\otimes v$.
\end{lemma}

\paragraph{A first example.} Let $1\le k\le d$ be an integer, and consider the function $M\longmapsto P_k(M):=\sigma_k({\rm Sp}\,M)$, where ${\rm Sp}\,M$ denotes the spectrum of the matrix $M$ and $\sigma_k$ is the elementary symmetric polynomial in $d$ variables. For instance $P_1$ is the trace and $P_d$ is the determinant. One has $P_k\ge0$ over the cone ${\bf Sym}_d^+$. If $k\ge2$, we have in addition $P_k(v\otimes v)=0$ for every vector $v$. Actually
\begin{equation}
\label{eq:Pksig}
P_k(m\,e\otimes e+(d-1-m)I_d)\sim(d-1)(d-1-m)^{k-1}.
\end{equation}

If $f_k:=P_k^{\frac d{k(d-1)}}$ exhibits a maximal gain of regularity, then (\ref{eq:Pksig}) together with Lemma \ref{l:gTm} yield
$$\frac{d(k-1)}{k(d-1)}\,\ge1,$$
that is $k\ge d$. Therefore only $f_d=\det^{\frac1{d-1}}$ exhibits a maximal gain of regularity among the class of functions $f_k$.

\subsection{The case of immanants}

The determinant is not the only homogeneous polynomial to vanish on the cone of rank-one tensors. Most of the immanants\footnote{There are two spellings in the litterature. {\em Immanant} recalls the determinant, while {\em immanent} recalls the permanent.} share this property. We recall their definition:
\begin{defin}
Let $G$ be a subgroup of the symmetric group ${\mathfrak S}_d$ and $\chi$ be an irreducible character over $G$. The immanant $J^G_\chi:{\bf M}_d(\R)\longrightarrow\R$ is the polynomial
$$J^G_\chi(M)=\sum_{g\in G}\chi(g)\prod_{i=1}^dm_{ig(i)}.$$
\end{defin}

If $G=\mathfrak S_d$, an immanant has integral coefficients. Two examples of immanants are the determinant, for which $\chi$ is the signature, and the {\em permanent}, denoted ${\rm Per}$, where $\chi\equiv1$. When $G=(1)$, there is only one immanant, the diagonal product.

When $G$ is a proper subgroup, and  if there exists a $g\in G$ such that $g^{-1}$ and $g$ are in different conjugacy classes, an immanant may have non-real coefficients.  However, the restriction of an immanant to ${\bf Sym}_d$ does have real coefficients~; the reason is that $\chi(g^{-1})=\overline{\chi(g)}$, while for symmetric matrices
$$\prod_{i=1}^ds_{ig^{-1}(i)}=\prod_{i=1}^ds_{ig(i)}.$$

In \cite{Sch}, I. Schur proved the following property of immanants over the cone of positive semi-definite matrices. We use the digit $1$ to denote the unit element (the identity) of $G$.
\begin{prop}[I. Schur, 1918.]
For every $S\in{\bf Sym}_d^+$, one has
$$\det S\le\frac1{\chi(1)}\,J^G_\chi(S).$$
\end{prop}

The particular case when $G=(1)$,
$$\det S\le s_{11}\cdots s_{dd},\qquad\forall\,S\in{\bf Sym}_d^+,$$
is known as the Hadamard inequality. The {\em permanent dominance conjecture} tells on the contrary that
$$\frac1{\chi(1)}\,J^G_\chi(S)\le{\rm Per}(S),\qquad\forall\,S\in{\bf Sym}_d^+.$$

Schur's property allows us to define 
$$f^G_\chi(S)=\left(J^G_\chi(S)\right)^{\frac1{d-1}},$$
which is homogeneous of degree $\frac d{d-1}$ over ${\bf Sym}_d^+$.  If $v\in\R^n$, we have
$$J^G_\chi(v\otimes v)=\left(\sum_{g\in G}\chi(g)\right)v_1^2\cdots v_d^2.$$
For $\chi\not\equiv1$,  the sum above equals zero by orthogonality with the trivial character. Therefore $f^G_\chi$ vanishes identically over the cone of tensors $v\otimes v$, and the question of whether $f^G_\chi$ displays a maximal gain of integrability becomes natural. 
The answer is given in Theorem \ref{th:imm} below. To begin with, we examine the degree of the polynomial $p^G_\chi(X):= J^G_\chi(I_d+X\,e\otimes e)$.
\begin{lemma}\label{l:pchi}
The degree of $p^G_\chi$ is $\ge2$, unless $G=\mathfrak S_n$ and $\chi=\epsilon$ is the signature.
\end{lemma}

Remark that this degree is $1$ for the pair $({\mathfrak S}_d,\epsilon)$.

\bigskip

\bepr

The coefficient of $X^2$ in $p^G_\chi$,
$$h^G_\chi(e):=\chi(1)\sum_{i<j}e_i^2e_j^2+\sum_{\tau=(i\,j)\in G}\chi(\tau)e_i^2e_j^2,$$
is a linear combination of monomials $e_i^2e_j^2$. The coefficient of a monomial is either $\chi(1)+\chi(\tau)$ or $\chi(1)$, depending on whether the transposition $\tau=(i\,j)$ belongs to $G$ or not (notice that $\chi(\tau)$ is real because $\tau^{-1}=\tau$).

Let $\rho:G\longrightarrow GL(V)$ denote the representation associated with $\chi$. We recall that  $|\chi(g)|\le\chi(1)$ for every $g\in G$, and the equality implies that $\rho(g)$ is a homothety,  $\rho(g)=\frac{\chi(g)}{\chi(1)}\,{\rm id}_V$. In particular $h^G_\chi$ has non-negative coefficients, from which we deduce $h^G_\chi\ge0$. 

If $\deg p^G_\chi\le1$, we have $h^G_\chi\equiv0$. Since $\chi(1)>0$, this implies that $G$ contains every transposition, and therefore $G=\mathfrak S_d$. It also tells us that $\chi(\tau)=-\chi(1)$ for every transposition, which implies $\rho(\tau)=-{\rm id}_V$. Because $\rho$ is a morphism and $\epsilon(g)=(-1)^\ell$ where $\ell$ is the number of transposition in a factorisation of $g$, we deduce $\rho(g)=\epsilon(g){\rm id}_V$. Hence $\chi=\epsilon$.

\enpr

\bigskip

This yields the main result of this section.
\begin{thm}\label{th:imm}
Among the pairs $(G,\chi)$, only one yields a function $f^G_\chi$ that displays a maximal gain of integrability, namely the pair $({\mathfrak S}_d,\epsilon)$, for which the function is $\det^{\frac1{d-1}}$.
\end{thm}

\bepr

Denote $h_k(e)X^k$ the dominant monomial in $p^G_\chi$. When $m\longrightarrow (d-1)^-$, we have
$$J^G_\chi(m\,e\otimes e+(d-1-m)I_d)\sim(d-1-m)^{d-k}m^kh_k(e).$$
In the formula
$$h_k(e)=\sum_{|A|=k}e_A^2\sum_{g\in G\,|\,A^c\subset{\rm Fix}(g)}\chi(g),\qquad e_A=\prod_{i\in A}e_i,$$
each sum of characters is  a non-negative integer. For the sake of completeness, it equals the multiplicity of the trivial representation in ${\rm Res}^G_H(\rho)$, where $H$ is the subgroup of all $g\in G$ such that $A^c\subset{\rm Fix}(g)$ (the stabilizer of $A^c$).

By definition $h_k\not\equiv0$ and we obtain $\int_{S^{d-1}}(h_k(e))^{\frac1{d-1}}ds(e)>0$. Therefore 
$$\int_{S^{d-1}}f^G_\chi(m\,e\otimes e+(d-1-m)I_d)\,ds(e)\sim\mu(d-1-m)^{\frac{d-k}{d-1}}$$
as $m\longrightarrow(d-1)^-$, where $\mu$ is some positive constant. If $f^G_\chi$ displays a maximal gain of integrability, this together with Lemma \ref{l:gTm} implies $\frac{d-k}{d-1}\ge1$, that is $k\le1$. 

We conclude by using Lemma \ref{l:pchi}.

\enpr

\section{More examples of divergence-free symmetric tensors}\label{s:otherDPT}

\subsection{The wave equation}

Let $u$ be a solution of the wave equation~: 
$$\partial_t^2u=c^2\Delta_y u$$
The assumption of bounded energy yields the standard regularity $\nabla_{t,y}u\in L^\infty(\R;L^2(\R^n))$.
Then the tensor
$$T=\begin{pmatrix}
\frac12(u_t^2+c^2|\nabla u|^2) & -c^2u_t\nabla u^T \\ -c^2u_t\nabla u &
c^4\nabla u\otimes\nabla u+\frac{c^2}2(u_t^2-c^2|\nabla u|^2)I_n
\end{pmatrix}.$$
is integrable in slabs $(0,\tau)\times\R^n$. We leave the reader checking that ${\rm Div}_{t,y}T=0$~; the first line of of this identity is nothing but the conservation of energy.

The matrix $T(t,y)$ admits two invariant subspaces, namely the plane $\R\times\R\nabla u$ and its orthogonal. The restriction of $T$ to $\R\times\R\nabla u$ is positive definite. In particular, $T$ is a DPT when $n=1$. On the contrary, the restriction to $\{0\}\times\nabla u^\bot$ is a homothety of ratio $\frac{c^2}2(u_t^2-c^2|\nabla u|^2)$. 
When $n\ge2$, this tensor is therefore positive semi-definite if and only if $c|\nabla u|\le|u_t|$.
One has
$$\det T=c^{2n}\left(\frac{u_t^2-c^2|\nabla u|^2}2\right)^{n+1}.$$
Notice that this determinant has a constant (positive) sign when the space dimension is odd, even if $T$ is not positive semi-definite.

\bigskip

If the compensated integrability held true, no matter of the signature of the tensor,  we should receive the estimate
$$\int_\R\int_{\R^n}\left|u_t^2-c^2|\nabla u|^2\right|^{1+\frac1n}dy\,dt\le CE_0^{1+\frac1n},$$
where $$E_0\equiv\int_{\R^n}\frac12(u_t^2+c^2|\nabla u|^2)\,dy$$ is the total energy. If $n=1$, this estimate is correct since $T$ is a DPT, and the proof is actually an application of Fubini's Theorem to Riemann invariants. On the contrary, in space dimension $n=2$ or $3$, Klainerman \& Machedon constructed an example where the energy is finite, but   $u_t^2-c^2|\nabla u|^2\not\in L^{1+\frac1n}_{t,y}$. See Proposition 5 in  \cite{KM}. Therefore the compensated integrability may fail when we drop the positivity assumption.

\paragraph{Remark.} If we consider the Laplace equation instead, $\Delta_xu=0$, the symmetric tensor $T=\nabla u\otimes\nabla u-\frac12\,|\nabla u|^2I_d$ is again divergence-free, and neither positive nor negative. But contrary to the wave equation,  the gain of integrability does hold true. Of course, this is not a manifestation of compensated integrability. The normal trace 
$$T\vec n=\frac{\partial u}{\partial n}\,\nabla u-\frac12\,|\nabla u|^2\vec n$$
satisfies $|T\vec n|=\frac14\,|\nabla u|^2$. The assumption that $T\vec n$ is a bounded measure over $\partial\Omega$ tells us that the normal derivative belongs to $L^2(\partial\Omega)$. By regularity theory for elliptic equations, we obtain $u\in H^{3/2}(\Omega)$. Then the Sobolev embedding provides $\nabla u\in L^p(\Omega)$ with $p=\frac{2d}{d-1}$. Therefore 
$$|\det T|^{\frac1{d-1}}=\left(\frac12\,|\nabla u|^2\right)^{\frac d{d-1}}\in L^1(\Omega).$$

\subsection{The nonlinear Maxwell system}

For the sake of simplicity, we choose a system of physical units in which the speed of light is $c=1$.

\bigskip

The variational formulation of the Maxwell's system for the electro-magnetic field in vacuum is
$$\delta{\cal L}=0, \qquad\hbox{ where}\quad{\cal L}:=\int\!\!\int L(B,E)\,dy\,dt.$$
The ambient space is that of closed $2$-forms $\omega:=E\cdot(dt\wedge dy)+B\cdot(dy\wedge dy)$ on the $1+3$-dimensional Minkowski space with coordinates $(t,y)$. The closedness reads
$$\partial_tB+{\rm curl}\,E=0,\qquad {\rm div}\,B=0$$
while the variational principle yields
$$\partial_t D-{\rm curl}\,H=0,\qquad{\rm div}\,D=0\qquad\hbox{where}\quad D:=\frac{\partial L}{\partial E}\,,\qquad H:=-\frac{\partial L}{\partial B}\,.$$
Assuming that $(B,E)\mapsto(B,D)$ is a change of variable, we may define the energy density
$$W(B,D):=D\cdot E-L(B,E).$$
Then we have
\begin{equation}
\label{eq:ColDil}
E=\frac{\partial W}{\partial D}\,,\qquad H:=\frac{\partial W}{\partial B}.
\end{equation}
The relations (\ref{eq:ColDil}) were the starting point of the thermodynamical approach for non-linear models by Coleman \& Dill \cite{CD}.
Their general form of the Maxwell's system,
$$\partial_tB+{\rm curl}\frac{\partial W}{\partial D}=0,\quad\partial_tD-{\rm curl}\frac{\partial W}{\partial B}=0,\quad{\rm div}\,B={\rm div}\,=0$$
implies  four extra conservation laws
\begin{eqnarray}
\label{ener}
\partial_tW+{\rm div}(E\times H) & = & 0, \\
\label{poyn}
\partial_t(D\times B)-{\rm Div}(E\otimes D+H\otimes B)+\nabla(D\cdot E+B\cdot H-W) & = & 0.
\end{eqnarray}
The first equation of the system above justifies the terminology that $Q:=E\times H$ is the {\em electromagnetic momentum}. The vector field $P:=D\times B$ is called the {\em Poynting vector}.

\bigskip

When the Maxwell's equations are consistent with special relativity, the Lagrangian $\cal L$ must be invariant under the action of the Lorentz group over $2$-forms. This means that the density $L$ depends only on two scalar variables $(\sigma,\pi)$,
$$L(B,E)=\ell(\sigma,\pi),\qquad\sigma:=\frac12(|B|^2-|E|^2),\quad\pi:=B\cdot E.$$
We have therefore
$$D=-\ell_\sigma E+\ell_\pi B,\qquad H=-\ell_\sigma B-\ell_\pi E$$
and this implies
\begin{equation}\label{PeqQ}
P=Q.
\end{equation}
This fundamental identity tells that the Abraham form $E\times H$ of the Poynting vector coincides with its Minkowski form $D\times B$.

\bigskip

Now, the system (\ref{ener},\ref{poyn}) can be rewritten in the form
$${\rm Div}_{t,y}T=0$$
where $T$ is a symmetric tensor, thanks to (\ref{PeqQ})~:
$$T=\begin{pmatrix}
W & P^T \\ P & \ell_\sigma(E\otimes E+B\otimes B)+(\ell+B\cdot H)I_3
\end{pmatrix}.$$

\paragraph{Decomposition of $T$.}

When $(B,E)$ are linearly independent, the $4\times4$ matrix $T$ admits two obvious invariant planes, namely $\R\times\R P$ and its orthogonal $\{0\}\times{\rm Span}(B,E)$. We therefore have $T\sim S\oplus^\bot R$, where
$$S=\begin{pmatrix} W & |P|^2 \\ 1 & \ell +B\cdot H \end{pmatrix}, \qquad R=\begin{pmatrix} \ell-\pi\ell_\pi & \pi\ell_\sigma \\ \pi\ell_\sigma & \ell-2\sigma\ell_\sigma-\pi\ell_\pi  \end{pmatrix}.$$
The spectrum of $T$, which is real, is therefore the union of that of $S$ and $R$. We have
$$\det S=-\det R=\ell_\sigma^2(\sigma^2+\pi^2)-(\ell-\sigma\ell_\sigma-\pi\ell_\pi)^2.$$
Thus $\det T=-(\det S)^2$ and $T$ cannot be positive semi-definite. It cannot be negative either because in practice $W$ is positive.

In conclusion, we cannot expect a compensated integrability for the tensor $T$.

\subsection{Systems of conservation laws with a convex extension}

The space-time dimension is still $d=1+n$, while the field $u$ takes values in a convex domain of $\R^N$. We consider a system of the form
$$\partial_tu+\sum_\alpha\partial_\alpha f^\alpha(u)=0$$
and assume that it is formally compatible with an additional conservation law
$$\partial_t\eta(u)+{\rm div}_x\vec Q(u)=0$$
where ${\rm D}^2\eta>0_N$ (strong convexity). Following Godunov \cite{God}, we introduce the conjugate variables $q_j=\partial\eta/\partial u_j$ and we know that there exist potentials $L^\alpha$ ($\alpha=0,\ldots,n$) such that $u=\partial L^0/\partial q$ and $f^\alpha(u)=\partial L^\alpha/\partial q$. The system rewrites
$$\partial_t\frac{\partial L^0}{\partial q}+\sum_\alpha\partial_\alpha\frac{\partial L^\alpha}{\partial q}=0.$$
The entropy and its fluxes are given by
$$\eta(u)=q\cdot u-L^0(q),\qquad Q^\alpha(u)=q\cdot\frac{\partial L^\alpha}{\partial q}-L^\alpha(q).$$

Let us assume $N=n+1=:d$. We label the indices $j$ from $0$ to $n$. Then ${\rm Div}_{t,y}T=0$ where $T=\partial L/\partial q$ is a Jacobian matrix. It is symmetric if and only if $L$ derives from a single potential $R(q)$:
$$L^\alpha=\frac{\partial R}{\partial q_\alpha}\,,\qquad T={\rm D}^2_qR.$$
It is positive if and only if $q\mapsto R(q)$ is a convex function. When it is so, Theorem \ref{th:slab} yields an estimate of
$$\int_0^\tau\int_{\R^n}(\det{\rm D}^2_qR)^{\frac1n}dx\,dt$$
in terms of the $L^1$-norm of the $\frac{\partial L^0}{\partial q}$ at initial and final times. In practical situations, these norms can be estimated in terms of the integral
$$I_0:=\int_{\R^n}\eta(u_0(y))\,dy$$
at initial time, using the second principle of thermodynamics, which tells that
$$t\mapsto I(t)=\int_{\R^n}\eta(u(t,y))\,dy$$
 is a non-increasing function of time. 
 
 \bigskip

The above analysis applies also when $N>d$ and the $d$ first lines of the system write ${\rm Div}_{t,y}T=0$ for a positive symmetric tensor. Then $R$ depends upon the additional parameters $(q_{d},\ldots,q_{N-1})$, which are not essential in the calculations.

\paragraph{Example: barotropic gas dynamics.}
In terms of the density $\rho$, the velocity $v$ and the linear momentum $m=\rho v$, we have
$$u=\begin{pmatrix} \rho \\ m \end{pmatrix},\qquad f^\alpha(u)=\begin{pmatrix} m_\alpha \\ m_\alpha v+p(\rho)\vec e^\alpha \end{pmatrix},\qquad \eta(u)=\frac{|m|^2}{2\rho}+\varepsilon(\rho).$$
This yields $q=(-\frac12\,|v|^2+\varepsilon'(\rho)\,,\,v)$ and
$$L^0(q)=\rho\varepsilon'-\varepsilon,\qquad L^\alpha=(\rho\varepsilon'-\varepsilon)v_\alpha.$$
Of course, we find that $L$ derives from a potential $R$, with
$$dR=(\rho\varepsilon'-\varepsilon)d(\varepsilon').$$
Let us denote $\lambda:=\varepsilon'(\rho)$. Then $\rho\varepsilon'-\varepsilon=\varepsilon^*(\lambda)$, where $\varepsilon^*$ denotes the convex conjugate of $\varepsilon$. We have $\rho=(\varepsilon^*)'(\lambda)$ and we find $R=R(\lambda)$ where $R'=\varepsilon^*$. With $\lambda=q_0+\frac12(q_1^2+\cdots+q_n^2)$, we have
$$T=R''V\otimes V+R'\begin{pmatrix} 0 & 0 \\ 0 &  I_n \end{pmatrix}, \qquad V:=\begin{pmatrix} 1 \\ q_1 \\ \vdots \\ q_n \end{pmatrix},$$
and the convexity of $R$ in terms of $q$ (the positivity of $T$) amounts to the convexity and monotonicity of $R$ as a function of $\lambda$.

To compute $\det T$, we observe that the plane ${\rm vec}(\vec e^0,V)$ and its orthogonal are $T$-invariant. On the former, $T$ acts as the matrix 
$$M=\begin{pmatrix} 0 & -R' \\ R'' & |V|^2R''+R' \end{pmatrix},$$
while on the latter $T$ acts as a homothety of ratio $R'$. Finally 
$$\det T=(R')^{n-1}\det M=(R')^nR''.$$
Therefore we obtain an estimate of
$$\int_0^\tau\int_{\R^n}R'(\lambda)R''(\lambda)^{\frac1n}dx\,dt.$$
We point out that $R'=\varepsilon^*=\rho\varepsilon'-\varepsilon$ is the pressure, while $R''=(\varepsilon^*)'=\rho$, since $(\varepsilon^*)'$ is the reciprocal of $\varepsilon'$.

\subsection{Donati compatibility condition}

In linearized elasticity, where typically $d=3$, the strain tensor $e$ associated with a displacement field $v$ is defined as the symmetric gradient,
$$e_{ij}=\frac12\,(\partial_iv_j+\partial_jv_i).$$
A natural question is to characterize the image of the map $v\longmapsto e$, which obviously depends upon the functional space chosen as the source. Besides the obvious symmetry $e_{ji}=e_{ij}$, a local compatibility condition due to Saint-Venant, 
\begin{equation}
\label{eq:SV}
\partial_i\partial_je_{k\ell}+\partial_k\partial_\ell e_{ij}=\partial_i\partial_\ell e_{jk}+\partial_j\partial_k e_{i\ell},
\end{equation}
is a collection of $d(d-1)\frac{d^2-4d+5}2$ differential identities.

When $\Omega$ is simply connected, the necessary conditions (\ref{eq:SV}) are also sufficient. For instance, if $e\in L^2(\Omega;{\bf Sym}_d)$ satisfies the Saint-Venant conditions, then there exists a vector field $v\in H^1(\Omega;\R^d)$ such that $e=\frac12(\nabla v+\nabla v^T)$ (see Ciarlet \& Ciarlet Jr \cite{CiCi} for the case $d=3$).

When $\Omega$ is a more general domain with an arbitrary topology, a sufficient condition must bear an integral form, in the spirit of De Rham Theorem for gradients fields. This integral form, called the {\em Donati compatibility condition}, writes
\begin{equation}
\label{eq:Donati}
\left(T:\Omega\stackrel{{\mathcal C}^\infty}\longrightarrow{\bf Sym}_d\,\hbox{ compactly supported }\,;\,{\rm Div}\,T=0\right)\Longrightarrow\left(\int_\Omega\Tr(Te)\,dx=0\right).
\end{equation}
For instance, Ting \cite{Ting} proved that if $e\in L^2(\Omega;{\bf Sym}_d)$ satisfies (\ref{eq:Donati}), then there exists a vector field $v\in H^1(\Omega;\R^d)$ such that $e=\frac12(\nabla v+\nabla v^T)$~; see also Ciarlet \& coll. \cite{CMM}. Boundary conditions can also be handled, see \cite{ACGK}.

\section{Vlasov models}\label{s:Vla}

Vlasov models are kinetic equations for clouds where particles don't collide, but are driven by a self-induced force field $F(t,y)$. In a plasma, the particles are ions, subject to the Coulomb force. In galaxies, the particles are stars and interact through gravity. The Newton equations
$$\frac{dX}{dt}=v,\qquad\frac{dv}{dt}=F$$
yield a kinetic model at the mesoscopic scale.
The density $f(t,y,v)$ obeys to the equation
$$(\partial_t+v\cdot\nabla_y)f+F(t,y)\cdot\nabla_vf=0.$$
The macroscopic density
$$\rho(t,y):=\int_{\R^n}f(t,y,v)\,dv$$
induces the force field:
$$F=-\nabla \phi,\qquad \phi:=\chi*\rho.$$
The radial kernel $\chi$ characterizes the physics at stake. The convolution is meant with respect to the space variable $y$,
$$\phi(t,y)=\int_{\R^n}\chi(|y-z|)\rho(t,z)\,dz=\int_{\R^n}\chi(|z|)\rho(t,y-z)\,dz.$$

In Coulomb or gravity force, the potential $\phi$ is given by a Poisson equation, $\Delta\phi={\rm cst}\cdot\rho$, and we speak of the Vlasov--Poisson model. Depending upon the sign of the constant, the force is repulsive or attractive.

As usual, integration with respect to the velocity variable yields mass conservation:
\begin{equation}\label{eq:Vmass}
\partial_t\rho +{\rm div}_ym=0,\qquad m:=\int_{\R^n}f(t,y,v)\,v\,dv.
\end{equation}
Integrating against $v\,dv$ gives formally
\begin{equation}\label{eq:Vmom}
\partial_t m+{\rm Div}_y\left(\int_{\R^n}fv\otimes v\right)-\rho F=0,
\end{equation}
which implies the conservation of the overall momentum, because of
$$\int_{\R^n}\rho F\,dy=-\int_{\R^n\times\R^n}\rho(y)\rho(z)\nabla\chi(y-z)\,dy\,dz,$$
and $\nabla\chi$ is odd.

\paragraph{The energy estimate.} We have formally
$$\partial_t\int_{\R^n}\frac{|v|^2}2\,f\,dv+{\rm div}_y\int_{\R^n}\frac{|v|^2}2\,fv\,dv=-\int_{\R^n}F\cdot\nabla_vf\,\frac{|v|^2}2\,dv=F\cdot m=-\nabla_y\phi\cdot m.$$
Therefore
$$\partial_t\int_{\R^n}\frac{|v|^2}2\,f\,dv+{\rm div}_y\left(\int_{\R^n}\frac{|v|^2}2\,fv\,dv+\phi m\right)=\phi\,{\rm div}_ym=-(\chi*\rho)\,\partial_t\rho.$$
Integrating in space, there comes
$$\frac d{dt}\,\int_{\R^n\times\R^n}\frac{|v|^2}2\,f\,dv\,dy=-\int_{\R^n}(\chi*\rho)\,\partial_t\rho\,dy.$$
Because $\chi$ is even, the bilinear from
$$(\rho,\eta)\longmapsto\int_{\R^n}(\chi*\rho)\,\eta\,dy$$
is symmetric. We infer
$$\frac d{dt}\,\int_{\R^n}(\chi*\rho)\,\rho\,dy=2\int_{\R^n}(\chi*\rho)\,\partial_t\rho\,dy,$$
and deduce formally
$$\frac d{dt}\,\left(\int_{\R^n\times\R^n}\frac{|v|^2}2\,f\,dv\,dy+\frac12\int_{\R^n}(\chi*\rho)\,\rho\,dy\right)=0.$$
The total energy
$$E(t):=\int_{\R^n\times\R^n}\frac{|v|^2}2\,f\,dv\,dy+\frac12\int_{\R^n}(\chi*\rho)\,\rho\,dy$$
is certainly a non-negative quantity if either $\chi$, or its Fourier transform $\hat\chi$ (which is a real function because $\chi$ is even), is non-negative. Actually, because of the conservation of mass, it suffices that $\chi$ be bounded by below ($\inf\chi(r)>-\infty$), in order to built a non-negative conserved quantity of the form $E+{\rm cst}\cdot M$. When this is the case, and the initial energy $E_0$ is finite, the conservation $E(t)\equiv E_0$ implies an {\em a priori} estimate.
In the sequel, we consider solutions that satisfy the conservation of mass and momentum (equations \ref{eq:Vmass} and \ref{eq:Vmom}), together with $E(t)\leq E_0$.

\paragraph{An alternate divergence-free tensor.}

Our next goal is two-fold. First we identify the expression $-\rho F$ as the space divergence of a symmetric tensor $S$, so that the following symmetric tensor be divergence-free,
$$T=\begin{pmatrix} \rho & m^T \\ m & \int_{\R^n}fv\otimes v\,dv+S \end{pmatrix}.$$
The second task is to determine those functions $\chi$ which guarantee that $S$ is non-negative. When this property occurs, $T$ is a DPT and we have the same estimate as for the Boltzmann equation, together with an estimate of
$$\int_0^\infty\int_{\R^n}(\rho\det S)^{\frac1n}dy\,dt$$
because of 
$$\det T=\rho\det\left(S+\int_{\R^n}fv\otimes v\,dv-\frac1\rho\,m\otimes m\right)\,\ge\rho\det S.$$
Notice that, because $\rho\mapsto F$ is linear, we expect that $S$ be quadratic in $\rho$. Therefore the new estimate involves a space-time integral of a quantity that is homogenenous in $\rho$, of degree $2+\frac1n$\,. This is a gain of integrability, compared to the energy estimate, which is only quadratic in $\rho$~; as usual, the price to pay is an extra integration, with respect to  the time variable. 

\bigskip

Dropping the time variable, we have
\begin{eqnarray*}
-(\rho F)(y) & = & (\rho\nabla \phi)(y)=\rho(y)\int_{\R^n}\nabla \chi(z)\rho(y-z)\,dz=\rho(y)\int_{\R^n}\chi'(|z|)\rho(y-z)\frac{z}{|z|}\,dz \\
& = & \frac12\,\rho(y)\int_{\R^n}\chi'(|z|)(\rho(y-z)-\rho(y+z))\frac{z}{|z|}\,dz
\end{eqnarray*}
At fixed $z$, we have
$$\rho(y)(\rho(y-z)-\rho(y+z))=-{\rm div}_y\int_{-\frac12}^{\frac12}\rho(y+(s-\frac12)z)\rho(y+(s+\frac12)z)z\,ds.$$
We infer $-\rho F={\rm Div}_yS$ with
$$S=-\frac12\,\int_{\R^n}\frac{\chi'(|z|)}{|z|}\,dz\int_{-\frac12}^{\frac12}\rho(y+(s-\frac12)z)\rho(y+(s+\frac12)z)z\otimes z\,ds.$$
In particular, $S$ is non-negative whenever $\chi$ is monotonous non-increasing. This corresponds to a repulsive force. This reminds us the case of the Euler equations, where the monotonicity of the pressure (the fact that the sound speed is a real number) induces a repulsive force too ; a gas tends to rarefy. Notice that the sign $\chi'\le0$ is also a necessary condition for the positivity~; think of a density $\rho$ that concentrates at two points.

\bigskip

According to Section \ref{s:slab}, a space-time estimate can be established whenever the DPT be integrable. For this we need that $S$ be integrable. This comes in general as a consequence of the mass and energy estimates. Actually, using the notation $\phi(r)=r|\chi'(r)|$, we have
$$\|S(y)\|\le\frac12\,\int_{\R^n}\phi(|z|)\,dz\int_{-\frac12}^{\frac12}\rho(y+(s-\frac12)z)\rho(y+(s+\frac12)z)\,ds.$$
This yields
$$\int_{\R^n}\|S(y)\|\,dy\le\frac12\,\int_{\R^n\times\R^n}\!\int_{-\frac12}^{\frac12}\phi(|z|)\,\rho(y+(s-\frac12)z)\rho(y+(s+\frac12)z)\,ds\,dy\,dz.$$
Let us make the change of variable
$$(y,z,s)\longmapsto(v,w,s),\qquad v=y+(s-\frac12)z,\,w=y+(s+\frac12)z.$$
Remarking that $dy\,dz\,ds=dv\,dw\,ds$, we obtain
$$\int_{\R^n}\|S(y)\|\,dy\le\frac12\,\int_{\R^n\times\R^n}\phi(|w-v|)\,\,\rho(v)\rho(w)\,dv\,dw.$$
Therefore $S\in L^\infty_t(L^1_y)$ whenever the initial data has finite mass and energy and the kernel satisfies
$$|\chi'(r)|\,r\le\hbox{ cst }\cdot(1+\chi(r)).$$

\bigskip

\paragraph{Comment.} The moral of this example is that we should investigate in more depth other natural examples of divergence-free symmetric tensors, when $T$ is neither positive nor negative semi-definite. We should ask ourselves whether there exists a hidden DPT $A$ such that the equation ${\rm Div}\,T=0$ is equivalent to ${\rm Div}\,A=0$. The latter tensor could be defined in terms of the state of the model in a non-local way, while $T$ is often defined in local terms.

\begin{english}

\end{english}

\end{document}